\documentclass[11pt]{amsart}

\usepackage{amssymb}
\usepackage{amsmath}
\usepackage{amsfonts}
\usepackage{amsthm}
\usepackage{amssymb, url}
\usepackage[mathscr]{eucal}
\usepackage{verbatim, color}
\usepackage{caption}
\usepackage{upgreek}

\usepackage[LY1]{fontenc}
\usepackage[utf8]{inputenc}
\usepackage{todonotes} 
\usepackage{pgf,tikz,pgfplots}
\pgfplotsset{compat=1.14}
\usepackage{mathrsfs}
\usetikzlibrary{arrows}

\newtheorem{theorem}{Theorem}[section]

\newtheorem{problem}[theorem]{Problem}
\newtheorem{remark}[theorem]{Remark}

\begin{document}
\title{On Yoshinaga's arrangement of lines and the containment problem}

\author{Maria Tombarkiewicz and Maciej Zi\k{e}ba}
\address{Department of Mathematics, Pedagogical University of Cracow,
	Podchor\c a\.zych 2,
	PL-30-084 Krak\'ow, Poland.}
\email{maria.tombarkiewicz@student.up.krakow.pl}
\address{Department of Mathematics, Pedagogical University of Cracow,
	Podchor\c a\.zych 2,
	PL-30-084 Krak\'ow, Poland.}
\email{maciej.zieba@student.up.krakow.pl}

\begin{abstract}
The main purpose of the note is to show that Yoshinaga's arrangement of $18$ lines having $48$ triple and $9$ double intersection points leads to a new (short) series of non-containment examples for  $I^{(3)} \subset I^{2}$, the question studied by Harbourne and Huneke.
\\
\textbf{Keywords}: line arrangements, symbolic powers of ideals, containment problem, intersection points \\
\textbf{AMS2010}: 13F20, 52C35, 32S22, 14N10
\end{abstract}

\maketitle

\section{Introduction}
In the recent years one can observe a lot of interest in comparing ordinary (or algebraic) and symbolic powers of homogeneous ideals. Let $\mathbb{K}$ be a field of characteristic $0$ and let $I \subset \mathbb{K}[x_{0}, ...,x_{n}]$ be a homogeneous ideal. The $r$-th algebraic power of $I$ is generated by $r$-fold products of all elements sitting in $I$, which is a purely algebraic concept. On the other side, we can look at the geometric side lurking behind the ideals, namely the concept of symbolic powers of  homogeneous ideals $I$. Assuming that $I$ is radical, by the celebrated Nagata-Zariski result we know that the $m$-th symbolic power can be interpreted as the set of all homogeneous forms in $n+1$ variables vanishing along ${\rm Zeros}(I)$ with multiplicities at least $m$. By definition, we see that $I^{m} \subseteq I^{(m)}$ for every $m\geq 1$, and it is natural to ask about the reverse inclusion.
\begin{problem}
Let $I \subset \mathbb{K}[x_{0}, ..., x_{n}]$ be a radical ideal, decide for which $m$ and $r$ there is the containment
$$I^{(m)} \subset I^{r}.$$
\end{problem}
The breakthrough has been achieved in the early 2000s when Ein, Lazarsfeld and Smith in characteristic zero \cite{ELS01} and Hochster and Huneke in positive characteristic \cite{HH} (see also Ma and Schwede \cite{MS} for mixed characteristic case) proved the following uniform relation.
\begin{theorem}
Let $I\subset \mathbb{K}[x_{0},...,x_{N}]$ be a homogeneous ideal such that every embedded component of its zero locus ${\rm Zeros}(I)$ has codimension at most $e$. Then the containment 
$$I^{(m)}\subset I^{r}$$ holds provided that  $m \geq er$.
\end{theorem}
In particular, if we restrict our attention to the case of $\mathbb{P}^{2}$, then the above result tells us that for a finite set of mutually distinct points $\mathcal{P} = \{P_{1}, ..., P_{s}\}$ and the associated radical ideal $I$ one always has 
$$I^{(2r)} \subset I^{r}.$$
In 2006, Huneke asked whether the uniform bound for the radical ideals associated with finite sets of mutually distinct points is tight.

\begin{problem}[Huneke]
Let $\mathcal{P} \subset \mathbb{P}^{2}$ be a finite set of mutually distinct points and $I$ the associated radical ideal. Does the containment $$(\star) : \quad I^{(3)} \subset I^{2}$$
hold?
\end{problem}
Huneke observed that $(\star)$ holds if $\mathbb{K}$ has characteristic $2$, but it was an open problem (almost for 8 years) whether the problem has an affirmative answer in general. It turned out, somehow surprisingly that the containment $(\star)$ \textbf{does not hold} in general. The first non-containment example was discovered by Dumnicki, Szemberg, and Tutaj-Gasi\'nska \cite{DST13} -- their example is based on the dual Hesse arrangement of $9$ lines and $12$ triple intersection points. Shortly afterwards a plethora of non-containment examples was revealed (see for instance \cite{BNC,dumek,Real}), and many of them come from the singular loci of certain (extreme in some sense) line arrangements in the complex projective plane. 

The main purpose of the present note is to add to the above list another non-containment example which is based on Yoshinaga's arrangement of $18$ lines \cite[Example 2.2]{Dimca} -- we shall describe the whole construction in the forthcoming section. As a small spoiler we can unveil the mystery standing behind Yoshinaga's construction -- this is an extremely interesting arrangement constructed via a clever deformation of the $6$-th Fermat arrangement (or CEVA arrangement) of $18$ lines which is given by the linear factors of the defining polynomial 
$$Q_{6}(x,y,z) = (x^{6}-y^{6})(y^{6}-z^{6})(z^{6}-x^{6}).$$
This Fermat arrangement has exactly $36$ triple and $3$ sixtuple intersection points. It is well-known that the whole family of Fermat arrangements, given by
$$Q_{n}(x,y,z) = (x^{n}-y^{n})(y^{n}-z^{n})(z^{n}-x^{n})$$
with $n\geq 3$ provides non-containment examples to $(\star)$, which is proved in \cite{HS}. This is an interesting phenomenon due to the fact that this is the only known infinite family of complex line arrangements without double points delivering non-containment examples.

The main result of this note can be formulated as follows.
\begin{theorem}
Let $\mathcal{P}$ be the singular locus of Yoshinaga's arrangement of $18$ lines and denote by $I$ the associated radical ideal. Then 
$$I^{(3)} \not\subseteq I^{2}.$$
Let $I_{3}$ be the radical ideal of the singular sublocus of Yoshinaga's arrangement consisting of only triple intersection points. Then still
$$I_{3}^{(3)} \not\subseteq I_{3}^{2}.$$
\end{theorem}
\begin{remark}
	It is worth pointing out that the above theorem shows an interesting phenomenon, namely the containment problem does not hold neither for the set of triple points nor for the set of double and triple points, which is very rare -- in most cases one must stick to the subset of triple intersection points of a given arrangement.
\end{remark}
\begin{remark}
In fact the non-containment holds for any set of points between the set of triple points and the set of all singular points of the arrangement.	
\end{remark}
\section{Yoshinaga's arrangement of lines}

We start with the $6$-th Fermat arrangement which is given by the following defining equation
\begin{equation}
Q(x,y,z) = (x^{6}-y^{6})(y^{6}-z^{6})(z^{6}-x^{6}).
\end{equation}
The picture below shows an idea standing behind the construction -- it cannot be realize over the real numbers due to the celebrated Sylvester-Gallai theorem.

\begin{figure}[h]
	\centering
	\includegraphics[width=2.4in]{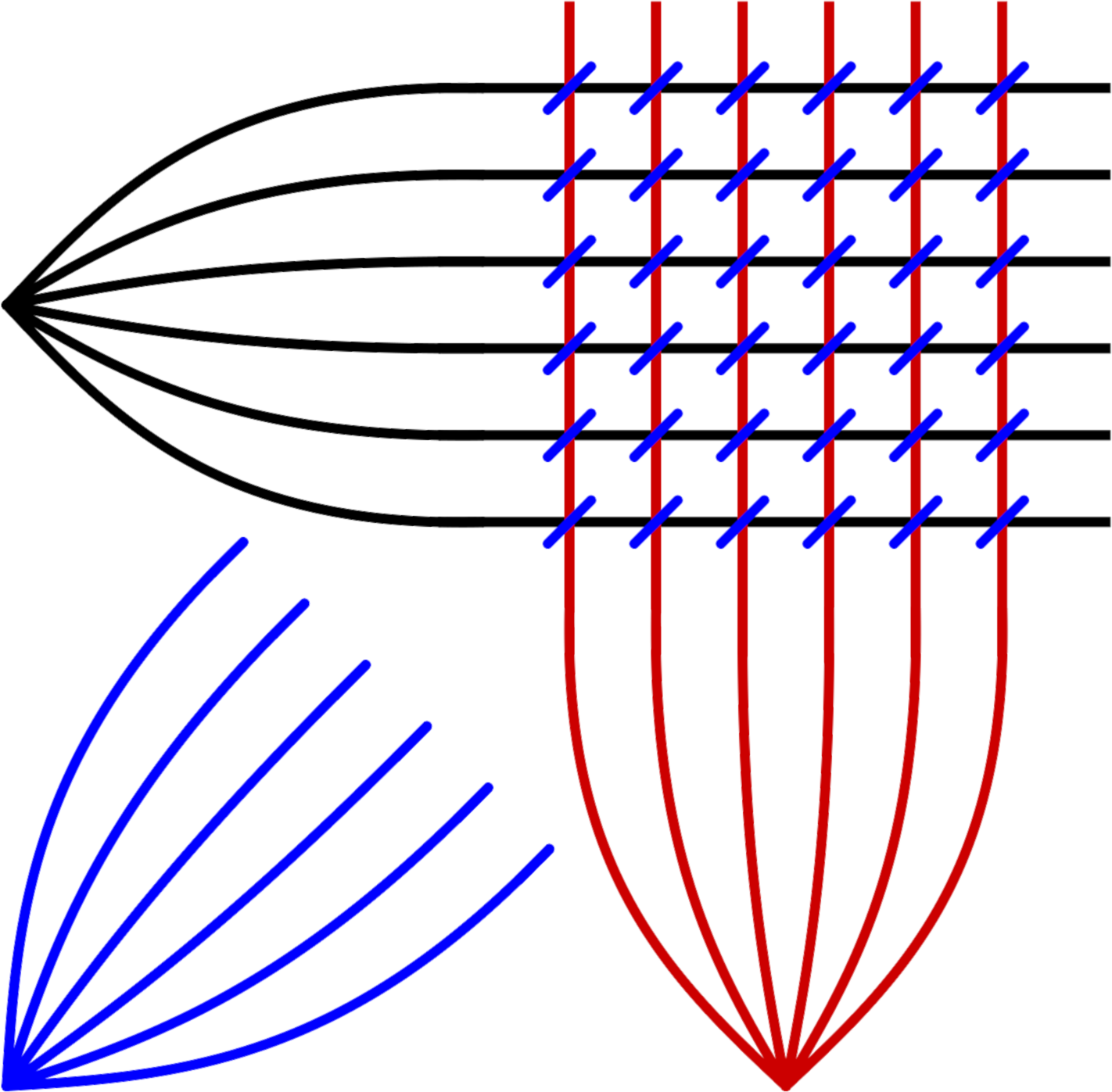}
	\caption{The $6$-th Fermat arrangement}
	\label{fig:figure1}
\end{figure}

Now we are going to present an outline of Yoshinaga's construction which is based on an interesting deformation argument.

%

Let $c\in \mathbb{R}$ be a large real number (for us it will be enough to take $c=15$). Set $a := e^{\frac{2\pi i}{6}}$. 
	Now we perform an appropriate deformation of three subpencils of lines, each consisting of $6$ lines intersecting simultaneous at a single point -- these points are the three reflectors in $6$-th Fermat arrangement. Moreover, this deformation is rather demanding due to the fact that we are going to maintain a \emph{complete-intersection-type} grid of $36$ triple points as in Fermat's construction.
	Our starting point is the polynomial which defines a sixtuple intersection point, let us take $P_1(x,y,z) = x^6-y^6$. Our deformation is given by the following polynomial
	\[
	P_{1}'(x,y,z)=(x^3-y^3)(x+y-cz)(ax+a^{5}y+cz)(a^{5}x+ay+cz),
	\]
	and it is easy to see that $P'$ defines an arrangement which can be viewed as the so-called $\mathcal{A}_{1}(6)$ simplicial arrangement of lines having exactly $6$ lines, $4$ triple points, and $3$ double points. Direct computations lead to 
	\[
	P_{1}'(x,y,z)=x^6-y^6+3cx^{4}yz-3cxy^{4}z-c^{3}x^{3}z^{3}+c^{3}y^{3}z^{3}
	\]
	
	Consider a cyclic permutation of variables $\uptau (x, y, z)=(y, z, x)$, then we define new polynomials in $\mathbb{C}[x,y,z]$ with respect to the action of $\uptau $- permutation, namely $P_{2}'(x,y,z) = \uptau P_{1}' = P_{1}'(y,z,x) $ and $P_{3}'(x,y,z) = \uptau^{2} P_{3}' = P_{1}'(z,x,y)$. 
	
	Set $P:=P_{1}'P_{2}'P_{3}'$. It turns out that $P$ splits into linear factors and delivers Yoshinaga's arrangement of lines. 
	
	\begin{figure}[ht]
		\centering
		\includegraphics[width=3in]{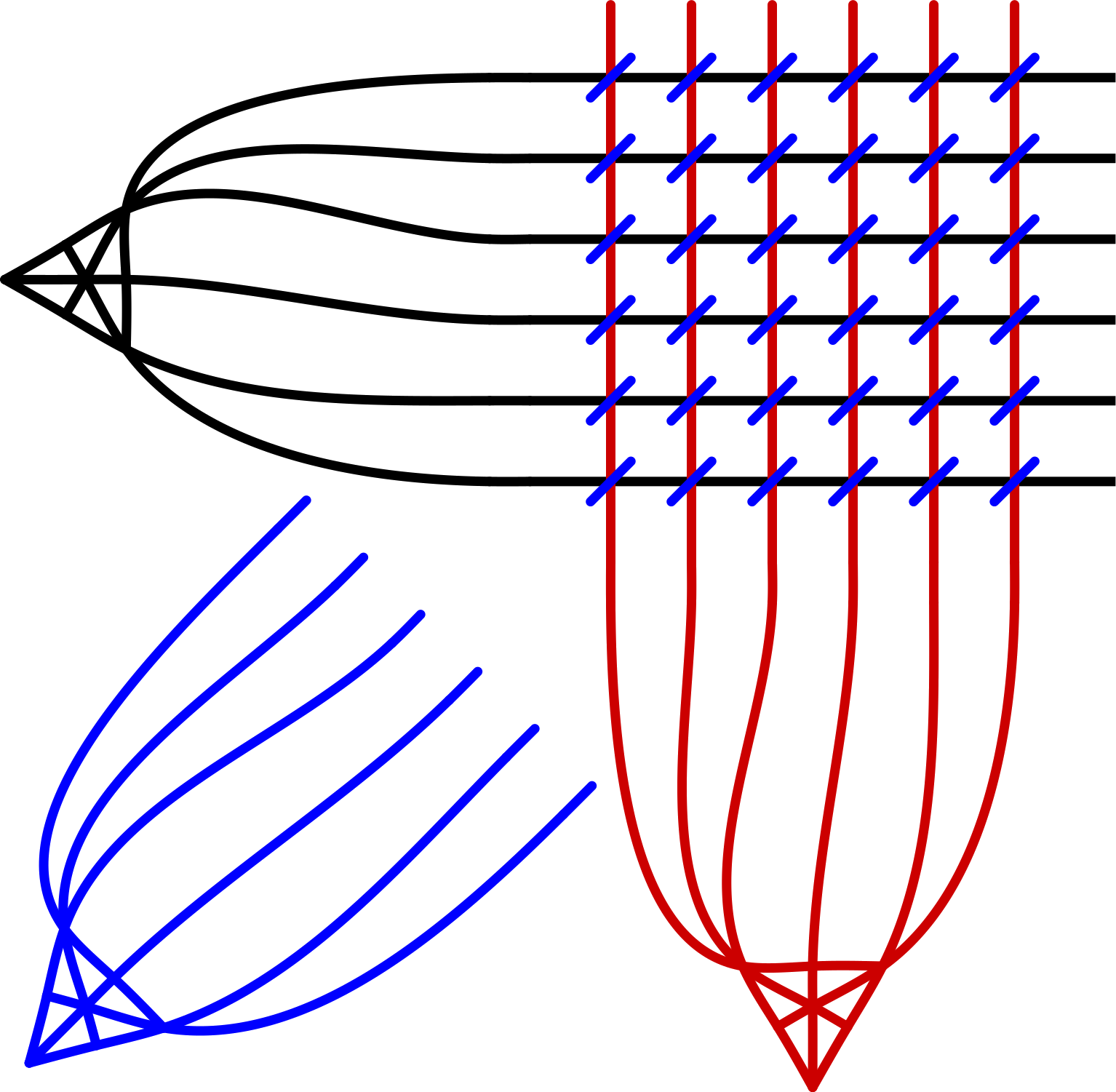}
		\caption{Yoshinaga's arrangement}
		\label{fig:figure2}
	\end{figure}

	We present below the equations of those $18$ lines (with $c=15$). 
\begin{align*} 
\ell_{1}&:  x-y & \ell_{10}&:  x+15ay+(a-1)z  \\
\ell_{2}&:  x+ay & \ell_{11}&:  x+(a-1)y+15az \\
\ell_{3}&:  x+(-a+1)y  &\ell_{12}&:  x-ay+(-15a+15)z \\
\ell_{4}&:  x-15y+z &\ell_{13}&:  x-z \\
\ell_{5}&:  x-\frac{1}{15}y-\frac{1}{15}z   &\ell_{14}&:  x+(-a+1)z \\   
\ell_{6}&:  x+y-15z & \ell_{15}&:  x+az  \\
\ell_{7}&:  x+\frac{1}{15}ay+\bigg(-\frac{1}{15}a+\frac{1}{15}\bigg)z &\ell_{16}&:  y-z \\ 
\ell_{8}&:  x+\bigg(-\frac{1}{15}a+\frac{1}{15}\bigg)y+\frac{1}{15}az &\ell_{17}&:  y+(-a+1)z \\
\ell_{9}&:  x+(-15a+15)y-az & \ell_{18}&:  y+az \\
\end{align*}
In Figure \ref{fig:figure2} we present the idea standing behind Yoshinaga's construction, please bear in mind that such an arrangement cannot be constructed over the reals.

\section{Non-containment}
Now we are going to describe the non-containment result. In the most cases we have the following yoga:
\begin{itemize}
	\item Consider an arrangement of lines $\mathcal{A} = \{\ell_{1}, ..., \ell_{d}\} \subset \mathbb{P}^{2}_{\mathbb{C}}$ having some extreme properties (in the sense of geometry or combinatorics), for us the crucial thing is that a given arrangement possesses large number of triple intersection points comparing with other singular points.
	\item Take the homogeneous form which is the product of defining equations of lines, i.e., if $\ell_{i} = V(f_{i})$, then $F = f_{1}\cdot ... \cdot f_{d}$. Show that $F \in I^{(3)} \setminus I^{2}$, where $I$ the radical ideal associated with a subset of singular points of $\mathcal{A}$ -- in most cases $I$ describes the set of triple and/or higher multiplicity intersection points.	
\end{itemize}

For the rest of this section, let us denote by $I$ the radical ideal of all $57$ singular points,  and by $I_{3}$ the radical ideal of all triple points in Yoshinaga's arrangement.

Here we followed a slightly different strategy, namely we compute explicitly generators of the second ordinary and the third symbolic power of a given ideal, reduce all generators of $I^{(3)}$ with respect to $I^{2}$, and we dig out a particular element which sits in $I^{(3)} \setminus I^{2}$. 
 A simple inspection of the aforementioned element reveals that it is the product (up to a non-zero constant) of the equations of $18$ lines defining Yoshinaga's arrangement plus a smooth curve of degree $3$ -- this curve is given by the following equation
\begin{equation}
\label{cubic}
x^3 + y^3 + z^{3} -\frac{3379}{225} xyz.
\end{equation}
This cubic curve is rather meaningful in that picture since it passes through all $9$ double intersection points with multiplicity $1$, but most importantly this is an element of the Hesse pencil with $(\mu:\lambda)=(1: -\frac{3379}{225})$ -- please consult \cite{AD} for details.
On the other side, it turns out that in the case of $I_{3}$ we can show that the element sitting in $I^{(3)}_{3} \setminus I^{2}_{3}$ is given exactly by the product of defining equations of $18$ lines, and this scenario fits into the picture that we have presented above. All computations are performed with use of Singular \cite{Singular} - please consult our script in the next section.
\begin{remark}
It is worth mentioning that a similar phenomenon was observed by Pokora and Ro\'e for Klein's arrangement of \,$21$ conic and \,$21$ lines in \cite{PR2019}. More precisely, the singular locus of the arrangement consists of $189$ quadruple points, $252$ triple points, and $42 = 21\cdot 2$ double intersection points -- these double points are exactly the intersections between pairs of lines and conics. If we denote by $I$ the radical ideal of all $483$ singular point, then a non-zero element from $I^{(3)} \setminus I^{2}$ giving the non-containment example is the product of the equations of $21$ lines, $21$ conics, and a smooth degree $6$ curve passing through $42$ double points -- this is the Hessian curve of Klein's quartic.
\end{remark}
\section{Singular Script}

{
\begin{verbatim}
LIB "primdec.lib";
option(redSB);
ring R=(0,e),(x,y,z),dp;
minpoly=e2-e+1;
proc rdideal(number pp, number qq, number rr) {
	matrix m[2][3]=pp,qq,rr,x,y,z;
	ideal I=minor(m,2);
	I=std(I);
	return(I);
}
list P2=15,15,2,15,(15e-15),(-2e),15,(-15e),(2e-2),2/15,1,1,(-2/15e),(e-1),1,
(2/15e-2/15),(-e),1,15,2,15,15,(2e-2),(-15e),15,(-2e),(15e-15);

list P3=0,0,(e+1),(15e-15),(15e-15),(-e+1),1/15,1/15,14/15,
(-1/15e),(-1/15e),(-1/15e+16/15),(1/15e-1/15),(1/15e-1/15),(1/15e+1),
1,1,1,(-e),(-e),1,(e-1),(e-1),1,(-15e),15,(-e+1),
(-1/15e),1/15,(-e-1/15),(1/15e-1/15),(-1/15e),(-16/15e+1/15),
1/15,(1/15e-1/15),(-14/15e),(-e),1,1,(e-1),(-e),1,1,(e-1),1,
(15e-15),15,(e),(1/15e-1/15),1/15,(e-16/15),1/15,(-1/15e),(14/15e-14/15),
(-1/15e),(1/15e-1/15),(16/15e-1),(e-1),1,1,1,(-e),1,(-e),(e-1),1,
-16/15,-224/15,-16/15,(-14/15e+1),(-1/15e-15),(-1/15e-14/15),
(14/15e+1/15),(1/15e-226/15),(1/15e-1),14,1,1, (e+15),(-e),1,
(-e+16),(e-1),1,(-14/15e+1),(-15e+226/15),(e-1/15),
(14/15e+1/15),(-226/15e+15),(14/15e+1/15),-16/15,(-224/15e+224/15),(16/15e),
(16e-15),1,1,(14e-14),(-e),1,(15e-16),(e-1),1,
(14/15e+1/15),(15e+1/15),(-e+14/15),-16/15,(224/15e),(-16/15e+16/15),
(-14/15e+1),(226/15e-1/15),(-14/15e+1),(-16e+1),1,1,(-15e-1),(-e),1,
(-14e),(e-1),1,(-2/225e+1/225),(-1/15e-1/15),(-1/15e+2/15),
(2/225e-1/225),(1/15e-2/15),(1/15e+1/15),(e+1),0,0,
(-2/225e+1/225),(2/15e-1/15),(2/15e-1/15),(-30e+15),(-e+2),(15e+15),
(30e-15),(e+1),(-15e+30),(-30e+15),(2e-1),(-30e+15),0,(-e-1),0;

int i;ideal I=1;ideal I3=1;
"";"generating ideals I^(3) and I^2 for triple points";"";
for(i=1;i<=(size(P3) div 3);i++){
	I=intersect(I,rdideal(P3[3*i-2],P3[3*i-1],P3[3*i]));
	I3=intersect(I3,rdideal(P3[3*i-2],P3[3*i-1],P3[3*i])^3);
}
I=std(I^2);I3=std(I3);
"";"reduction";"";
NF(I3,I);
"";"decomposition of an element from I^(3)\setminus I^2";"";
factorize(I3[1]);

I=1;I3=1;
"";"generating ideals I^(3) and I^2 for double and triple points";"";
P3=P3+P2;
for(i=1;i<=(size(P3) div 3);i++){
	I=intersect(I,rdideal(P3[3*i-2],P3[3*i-1],P3[3*i]));
	I3=intersect(I3,rdideal(P3[3*i-2],P3[3*i-1],P3[3*i])^3);
}
I=std(I^2);I3=std(I3);
"";"reduction";"";
NF(I3,I);
"";"decomposition of an element from I^(3)\setminus I^2";"";
factorize(I3[1]);
\end{verbatim}}

\section*{Acknowledgments}
We would like to warmly thank Prof. Masahiko Yoshinaga for explaining his construction and Grzegorz Malara for his help with Singular computations. The second author is partially supported by the Polish Ministry of Science and Higher Education within the program ``Najlepsi z najlepszych 4.0'' (''The best of the best 4.0'').\\
\newpage
\section*{Appendix}
Here we present the table of incidences between lines which gives the combinatorics of triple intersection points. 

\begin{table*}[hb]
	\scalebox{0.7}
	{
\centering
\begin{tabular}{|l|l|l|l|l|l|l|l|l|l|l|l|l|l|l|l|l|l|l|} 
\hline
$\cdot$ & $\ell_{1}$ & $\ell_{2}$ & $\ell_{3}$ & $\ell_{4}$ & $\ell_{5}$ & $\ell_{6}$ & $\ell_{7}$ & $\ell_{8}$ & $\ell_{9}$ & $\ell_{10}$ & $\ell_{11}$ & $\ell_{12}$ & $\ell_{13}$ & $\ell_{14}$ & $\ell_{15}$ & $\ell_{16}$ & $\ell_{17}$ & $\ell_{18}$  \\ 
\hline
$P_{1}$      & +  & +  & +  &    &    &    &    &    &    &     &     &     &     &     &     &     &     &      \\ 
\hline
$P_{1}$      & +  &    &    &    & +  & +  &    &    &    &     &     &     &     &     &     &     &     &      \\ 
\hline
$P_{2}$      & +  &    &    &    &    &    & +  &    &    &     &     &     & +   &     &     &     &     &      \\ 
\hline
$P_{3}$      & +  &    &    &    &    &    &    & +  &    &     &     &     &     & +   &     &     &     &      \\ 
\hline
$P_{4}$      & +  &    &    &    &    &    &    &    & +  &     &     &     &     &     & +   &     &     &      \\ 
\hline
$P_{5}$      & +  &    &    &    &    &    &    &    &    & +   &     &     &     &     &     & +   &     &      \\ 
\hline
$P_{6}$      & +  &    &    &    &    &    &    &    &    &     & +   &     &     &     &     &     & +   &      \\ 
\hline
$P_{7}$      & +  &    &    &    &    &    &    &    &    &     &     & +   &     &     &     &     &     & +    \\ 
\hline
$P_{8}$      &    & +  &    & +  &    & +  &    &    &    &     &     &     &     &     &     &     &     &      \\ 
\hline
$P_{9}$     &    & +  &    &    &    &    & +  &    &    &     &     &     &     &     & +   &     &     &      \\ 
\hline
$P_{10}$     &    & +  &    &    &    &    &    & +  &    &     &     &     & +   &     &     &     &     &      \\ 
\hline
$P_{11}$     &    & +  &    &    &    &    &    &    & +  &     &     &     &     & +   &     &     &     &      \\ 
\hline
$P_{12}$     &    & +  &    &    &    &    &    &    &    & +   &     &     &     &     &     &     & +   &      \\ 
\hline
$P_{13}$     &    & +  &    &    &    &    &    &    &    &     & +   &     &     &     &     &     &     & +    \\ 
\hline
$P_{14}$     &    & +  &    &    &    &    &    &    &    &     &     & +   &     &     &     & +   &     &      \\ 
\hline
$P_{15}$     &    &    & +  & +  & +  &    &    &    &    &     &     &     &     &     &     &     &     &      \\ 
\hline
$P_{16}$     &    &    & +  &    &    &    & +  &    &    &     &     &     &     & +   &     &     &     &      \\ 
\hline
$P_{17}$     &    &    & +  &    &    &    &    & +  &    &     &     &     &     &     & +   &     &     &      \\ 
\hline
$P_{18}$     &    &    & +  &    &    &    &    &    & +  &     &     &     & +   &     &     &     &     &      \\ 
\hline
$P_{19}$     &    &    & +  &    &    &    &    &    &    & +   &     &     &     &     &     &     &     & +    \\ 
\hline
$P_{20}$     &    &    & +  &    &    &    &    &    &    &     & +   &     &     &     &     & +   &     &      \\ 
\hline
$P_{21}$     &    &    & +  &    &    &    &    &    &    &     &     & +   &     &     &     &     & +   &      \\ 
\hline
$P_{22}$     &    &    &    & +  &    &    & +  &    &    &     &     &     &     &     &     & +   &     &      \\ 
\hline
$P_{23}$     &    &    &    & +  &    &    &    & +  &    &     &     &     &     &     &     &     &     & +    \\ 
\hline
$P_{24}$     &    &    &    & +  &    &    &    &    & +  &     &     &     &     &     &     &     & +   &      \\ 
\hline
$P_{25}$     &    &    &    & +  &    &    &    &    &    & +   &     &     & +   &     &     &     &     &      \\ 
\hline
$P_{26}$     &    &    &    & +  &    &    &    &    &    &     & +   &     &     &     & +   &     &     &      \\ 
\hline
$P_{27}$     &    &    &    & +  &    &    &    &    &    &     &     & +   &     & +   &     &     &     &      \\ 
\hline
$P_{28}$     &    &    &    &    & +  &    & +  &    &    &     &     &     &     &     &     &     & +   &      \\ 
\hline
$P_{29}$     &    &    &    &    & +  &    &    & +  &    &     &     &     &     &     &     & +   &     &      \\ 
\hline
$P_{30}$     &    &    &    &    & +  &    &    &    & +  &     &     &     &     &     &     &     &     & +    \\ 
\hline
$P_{32}$     &    &    &    &    & +  &    &    &    &    & +   &     &     &     &     & +   &     &     &      \\ 
\hline
$P_{33}$     &    &    &    &    & +  &    &    &    &    &     & +   &     &     & +   &     &     &     &      \\ 
\hline
$P_{34}$     &    &    &    &    & +  &    &    &    &    &     &     & +   & +   &     &     &     &     &      \\ 
\hline
$P_{35}$     &    &    &    &    &    & +  & +  &    &    &     &     &     &     &     &     &     &     & +    \\ 
\hline
$P_{36}$     &    &    &    &    &    & +  &    & +  &    &     &     &     &     &     &     &     & +   &      \\ 
\hline
$P_{37}$     &    &    &    &    &    & +  &    &    & +  &     &     &     &     &     &     & +   &     &      \\ 
\hline
$P_{38}$     &    &    &    &    &    & +  &    &    &    & +   &     &     &     & +   &     &     &     &      \\ 
\hline
$P_{39}$     &    &    &    &    &    & +  &    &    &    &     & +   &     & +   &     &     &     &     &      \\ 
\hline
$P_{40}$     &    &    &    &    &    & +  &    &    &    &     &     & +   &     &     & +   &     &     &      \\ 
\hline
$P_{41}$     &    &    &    &    &    &    & +  & +  &    &     & +   &     &     &     &     &     &     &      \\ 
\hline
$P_{42}$     &    &    &    &    &    &    & +  &    & +  &     &     & +   &     &     &     &     &     &      \\ 
\hline
$P_{43}$     &    &    &    &    &    &    &    & +  & +  & +   &     &     &     &     &     &     &     &      \\ 
\hline
$P_{44}$     &    &    &    &    &    &    &    &    &    & +   & +   & +   &     &     &     &     &     &      \\ 
\hline
$P_{45}$     &    &    &    &    &    &    &    &    &    &     &     &     & +   & +   &     &     & +   &      \\ 
\hline
$P_{46}$     &    &    &    &    &    &    &    &    &    &     &     &     & +   &     & +   &     &     & +    \\ 
\hline
$P_{47}$     &    &    &    &    &    &    &    &    &    &     &     &     &     & +   & +   & +   &     &      \\ 
\hline
$P_{48}$     &    &    &    &    &    &    &    &    &    &     &     &     &     &     &     & +   & +   & +    \\
\hline
\end{tabular}

}
		\caption*{Incidences of all triple points in Yoshinaga's arrangement}
\end{table*}


\newpage
\begin{thebibliography}{000}
\bibitem{AD}

Artebani, M., Dolgachev, I.:
The Hesse pencil of plane cubic curves. \textit{Enseign. Math.} (2) 55, No. 3-4: 235--273 (2009).

\bibitem{BNC}
Bauer, Th, Di Rocco, S., Harbourne, B., Huizenga, J., Lundman, A., Pokora, P., Szemberg, T.: Bounded negativity and arrangements of lines. \textit{Int. Math. Res. Notices} 2015: 9456 -- 9471 (2015).

\bibitem{Singular}
Decker, W.; Greuel, G.-M.; Pfister, G.; Sch{\"o}nemann, H.: 
{\sc Singular} {4-1-2} --- {A} computer algebra system for polynomial computations.
\url{http://www.singular.uni-kl.de} (2019).
\bibitem{Dimca}
A. Dimca, Monodromy of triple point line arrangements. Advanced Studies in Pure Mathematics 66, 2015. Singularities in Geometry and Topology 2011 pp. 71-80
\bibitem{dumek}
Dumnicki, M., Harbourne, B., Nagel, U., Seceleanu, A., Szemberg, T., Tutaj-Gasi\'nska, H.: Resurgences for ideals of special point configurations in $\mathbb{P}^{N}$ coming from hyperplane arrangements. \textit{J. Algebra} 443: 383--394 (2015).

\bibitem{Real}
Czapli\'nski, A., G\l \'owka, A., Malara, G., Lampa-Baczynska, M., \L uszcz-\'Swidecka, P., Pokora, P. and Szpond, J.:
A counterexample to the containment $I^{(3)}\subset I^2$ over the reals. \textit{Adv. Geom.} 16: 77--82 (2016).


\bibitem{DST13}
Dumnicki, M., Szemberg, T., Tutaj-Gasi\'nska, H.:
Counterexamples to the $I^{(3)}\subset I^2$ containment.
\textit{J.~Alg.} 393: 24--29 (2013).

\bibitem{ELS01}
Ein, L., Lazarsfeld, R., Smith, K.:
Uniform bounds and symbolic powers on smooth varieties.
\textit{Invent. Math}. 144: 241--252 (2001).

\bibitem{HS}
Harbourne, B., Seceleanu, A.: Containment counterexamples for ideals of various configurations of points in $\mathbb{P}^{N}$. \textit{J. Pure Appl. Algebra} 219: 1062--1072 (2015).

\bibitem{HH}
Hochster, M., Huneke, C.: Comparison of symbolic and ordinary powers of ideals. \textit{Invent. Math.} 147: 349--369 (2002).
\bibitem{MS}
Ma, L., Schwede, K.:  Perfectoid  multiplier/test  ideals  in  regular  rings  and  bounds  onsymbolic powers. \textit{Invent. Math.}  214(2): 913--955 (2018).
\bibitem{PR2019}
Pokora, P., Ro\'e. J.: The 21 reducible polars of Klein's quartic. \textit{Exp. Math.}, \url{https://doi.org/10.1080/10586458.2018.1488155}.
\end{thebibliography}
\end{document}